\documentclass{article}

\usepackage{amsfonts}
\usepackage{amssymb}
\usepackage{enumerate}
\usepackage[T1]{fontenc}

\newtheorem{ttt}{Theorem}[section]
\newtheorem{llll}[ttt]{Lemma}
\newtheorem{ccc}[ttt]{Claim}
\newtheorem{eee}[ttt]{Example}
\newtheorem{sss}[ttt]{Statement}
\newtheorem{ddd}[ttt]{Definition}
\newtheorem{qqq}[ttt]{Question}
\newtheorem{cccc}[ttt]{Corollary}

\newcommand{\bt}{\begin{ttt}}
\newcommand{\bl}{\begin{llll}}
\newcommand{\bc}{\begin{ccc}}
\newcommand{\bex}{\begin{eee}}
\newcommand{\bs}{\begin{sss}}
\newcommand{\bd}{\begin{ddd} \upshape}
\newcommand{\bq}{\begin{qqq}}
\newcommand{\bcor}{\begin{cccc}}

\newcommand{\bp}{\noindent\textbf{Proof. }}
\newcommand{\br}{\noindent\textbf{Remark. }}

\newcommand{\et}{\end{ttt}}
\newcommand{\el}{\end{llll}}
\newcommand{\ec}{\end{ccc}}
\newcommand{\eex}{\end{eee}}
\newcommand{\es}{\end{sss}}
\newcommand{\ed}{\end{ddd}}
\newcommand{\eq}{\end{qqq}}
\newcommand{\ecor}{\end{cccc}}

\newcommand{\ep}{\hspace{\stretch{1}}$\square$\medskip}

\newcommand{\RR}{\mathbb{R}}

\title{Less than $2^{\omega}$ many translates of a compact nullset may cover the real line}

\author{M\'arton Elekes}

\begin{document}

\maketitle 

\begin{abstract}
We answer a question of Darji and Keleti by proving in $ZFC$ that there exists a compact nullset $C_0\subset\RR$ such that for every perfect set $P\subset\RR$ there exists $x\in\RR$ such that $(C_0+x)\cap P$ is uncountable. Using this $C_0$ we answer a question of Gruenhage by showing that it is consistent with $ZFC$ that less than $2^\omega$ many translates of a compact nullset cover $\RR$.
\end{abstract}

\insert\footins{\footnotesize{MSC codes: Primary ; Secondary }}
\insert\footins{\footnotesize{Key Words: compact, measure zero, translate, cover, perfect, consistent}}



\bigskip
\bigskip

The results of this note will appear in a forthcoming paper with detailed proofs and some background. The main goal now is to state the results.

\bigskip
\bigskip

The following set is fairly well known, it was investigated for example by Erd\H os and Kakutani.

\bd Denote
\[
C_0 = \{\sum_{n=2}^{\infty} \frac{d_n}{n!} | d_n \in\{0,1,\dots,n-2\}\ \forall n \}
\]
\ed

It is easy to see that $C_0\subset\RR$ is a compact nullset, that is, a compact set of Lebesgue measure zero.
The following theorem answers a question of Darji and Keleti.

\bt
For every perfect set $P\subset\RR$ there exists a translate $C_0+x$ of $C_0$ such that $(C_0+x)\cap P$ is uncountable.
\et

\bp (sketch) Think of $d_n$ as digits with "increasing base", then every $x\in[0,1]$ has an (almost) unique expansion. Construct recursively a dyadic tree of intervals that describes a very thin perfect subset $Q\subset P$ which is "compatible" with this digit representation. Make sure that each branching is at a far enough level in order to be able to construct the expansion of $x$ such that $Q-x\subset C_0$ by simply counting the number of bad digits for every fixed $n$.
\ep

Now we can answer the original question of Gruenhage, which motivated the work of Darji and Keleti. $CPA$ is the well known axiom of Ciesielski and Pawlikowski.

\bt
CPA implies that $\RR$ can be covered by less than continuum many translates of $C_0$.
\et

\br (J. Steprans) We can also work with the so callad "slaloms" to obtain the same result for more models.

\bigskip

\noindent
\textsc{R\'enyi Alfr\'ed Institute,
Re\'altanoda u. 13-15. Budapest 1053, Hungary}

\textit{Email address}: \verb+emarci@renyi.hu+


\end{document}